\documentclass[11pt,twoside,reqno]{amsart}
\usepackage{mathptmx,amsmath,amssymb,amsfonts,amsthm,
mathptmx,enumerate,color}
\usepackage{graphicx}

\usepackage{epstopdf}

\newtheorem{theorem}{Theorem}

\newtheorem{lemma}{Lemma}

\theoremstyle{definition}

\newtheorem{remark}{Remark}

\numberwithin{equation}{section}

\newcommand{\qbinomial}[3]{\mbox{$
\biggl[ 
\begin{array}{c}
#1\\
 #2
\end{array}\biggr]_{
\!{#3}}$}}

\newcommand{\be}{\begin{equation}}
\newcommand{\ee}{\end{equation}}
\newcommand{\bea}{\begin{eqnarray}}
\newcommand{\eea}{\end{eqnarray}}
\newcommand{\bd}{\begin{displaymath}}
\newcommand{\ed}{\end{displaymath}}

\begin{document}

\setcounter{page}{1}

\vspace*{1.0cm}

\title[Generating Functions for Some Families of 
the Generalized Al-Salam-Carlitz $q$-Polynomials]{
Generating Functions for Some Families of the 
Generalized \\[1mm] Al-Salam-Carlitz $q$-Polynomials}

\author[H. M. Srivastava and Sama Arjika]
{H. M. Srivastava$^{1,2,3}$ and Sama Arjika$^{4,\ast}$}

\maketitle

\vspace*{-0.6cm}

\begin{center}
{\footnotesize {\it
$^{1}$Department of Mathematics and Statistics, 
University of Victoria, Victoria, British Columbia $V8W \; 3R4$, Canada\\
$^{2}$Department of Medical Research, 
China Medical University Hospital, China Medical University,\break
Taichung $40402$, Taiwan, Republic of China\\
$^{3}$Department of Mathematics and Informatics, Azerbaijan University,\\ 
71 Jeyhun Hajibeyli Street, AZ$1007$ Baku, Azerbaijan\\
$^{4}$Department of Mathematics and Informatics, University of Agadez,\\ 
Post Office Box $199$, Agadez, Niger}}
\end{center}

\vskip 4mm 

\begin{center} 
{\bf Abstract}
\end{center}
\begin{quotation}
In this paper, by making use of the familiar $q$-difference operators
$D_q$ and $D_{q^{-1}}$, we first introduce two homogeneous $q$-difference   
operators $\mathbb{T}({\bf a},{\bf b},cD_q)$  and  
$\mathbb{E}({\bf a},{\bf b}, cD_{q^{-1}})$, which turn out to be 
suitable for dealing with the families of the generalized Al-Salam-Carlitz 
$q$-polynomials $\phi_n^{({\bf a},{\bf b})}(x,y|q)$ and  
$\psi_n^{({\bf a},{\bf b})}(x,y|q)$. We then apply each of these
two homogeneous $q$-difference   
operators in order to derive generating functions, Rogers type formulas, the 
extended Rogers type formulas and the Srivastava-Agarwal type 
linear as well as bilinear generating 
functions involving each of these families of 
the generalized Al-Salam-Carlitz 
$q$-polynomials. We also show how the various results presented here
are related to those in many earlier works on the topics which we
study in this paper.
\end{quotation}   

\vskip 2mm


\vskip 1mm


\renewcommand{\thefootnote}{}
\footnotetext{$^{\ast}$Corresponding author.
\par
E-Mail: harimsri@math.uvic.ca (H. M. Srivastava);\; 
rjksama2008@gmail.com (Sama Arjika).} 

\section{Introduction, Definitions and Preliminaries}
The quantum (or $q$-) polynomials constitute a very interesting set of  
special functions and orthogonal polynomials. Their generating functions appear 
in several branches of mathematics and physics 
(see, for details, \cite{George86,GasparRahman,GMS-2001,JMedem,SrivastavaManocha}),
such as (for example) continued fractions, Eulerian series, theta functions, 
elliptic functions, quantum groups and algebras, discrete mathematics (including 
combinatorics and graph theory), coding theory, and so on.

Recently, new classes of special functions including (for example) 
$q$-hybrid special polynomials, $q$-Sheffer-Appell polynomials, twice-iterated 
2D $q$-Appell polynomials and a unified class of Apostol-type $q$-polynomials 
were introduced in \cite{YasminM2019},
\cite{YasminMu2021}, \cite{HariG2019} and \cite{HMS-BIMS2020} in which
some properties of the introduced polynomials were derived.  
For more information, the interested reader should 
refer to \cite{YasminM2019}, \cite{YasminMu2021}, \cite{HariG2019}
and \cite{HMS-BIMS2020}. 

In the year 1997, Chen and Liu \cite{Liu97}  
developed a method of deriving basic (or $q$-) hypergeometric identities 
by parameter augmentation, which may be viewed as being analogous to the 
method used rather extensively in the theory of {\it ordinary} 
hypergeometric functions and hypergeometric 
generating functions (see, for details, \cite{SrivastavaManocha}).  
Subsequent investigations along the lines developed in \cite{Liu97}
can be found in \cite{Abdlhusein2014}, \cite{Liu98}, \cite{Chen2008},
\cite{Chen2007}, \cite{Saad2014} and \cite{Srivastava-Arjika-Kelil}.  
The main object of this paper is to investigate two families of the generalized 
Al-Salam-Carlitz $q$-polynomials $\phi_n^{({\bf a},{\bf b})}(x,y|q)$ and  
$\psi_n^{({\bf a},{\bf b})}(x,y|q)$ by first representing them by 
the homogeneous $q$-difference operators $\mathbb{T}({\bf a},{\bf b},cD_q)$  
and $\mathbb{E}({\bf a},{\bf b},cD_{q^{-1}})$, which we have introduced here. 
We then derive a number of $q$-identities such as (among other results) 
generating functions, Rogers type formulas, two kind of the extended Rogers type formulas, 
and Srivastava-Agarwal type generating functions for each of the generalized 
Al-Salam-Carlitz $q$-polynomials $\phi_n^{({\bf a},{\bf b})}(x,y|q)$ and  
$\psi_n^{({\bf a},{\bf b})}(x,y|q)$.  

Here, in this paper, we adopt the common conventions and notations 
on $q$-series and $q$-hypergeometric functions. For the convenience of the 
reader, we provide a summary of the mathematical notations and definitions, 
basic properties and other relations to be used in the sequel.   
We refer, for details, to the general references 
(see \cite{Koekock,GasparRahman,SrivastavaKarlsson}) for the 
definitions and notations. Throughout this paper, 
we assume that  $|q|< 1$. 

For complex numbers $a$, the $q$-shifted factorials are defined by
\bea
(a;q)_0:=1,\quad  (a;q)_{n} =\prod_{k=0}^{n-1} (1-aq^k)   
\quad \text{and} \quad (a;q)_{\infty}:=\prod_{k=0}^{\infty}(1-aq^{k}),
\eea
where (see, for example, \cite{GasparRahman} and \cite{Slater})
$$(a;q)_n=\frac{(a;q)_\infty}{(aq^n;q)_\infty},\qquad 
(a;q)_{n+m}=(a;q)_n(aq^n;q)_m$$
and 
$$(aq^{-n};q)_n=(q/a;q)_n(-a)^nq^{-n-({}^n_2)}.$$

We adopt the following notation:  
$$(a_1,a_2, \cdots, a_r;q)_m=(a_1;q)_m (a_2;q)_m\cdots(a_r;q)_m
\qquad (m\in \mathbb{N}:=\{1,2,3,\cdots\}).$$
Also, for $m$ large, we have
$$
(a_1,a_2, \cdots, a_r;q)_\infty=(a_1;q)_\infty 
(a_2;q)_\infty\cdots(a_r;q)_\infty.$$
  
The $q$-numbers and the
$q$-factorials are defined as follows:
\be 
[n]_{q}:=\frac{1-q^n}{1-q},\quad [n]_{q}!
:=\prod_{k=1}^n\left(\frac{1-q^k}{1-q}\right)
\qquad \text{and} \qquad [0]_{q}!:=1.
\ee
The $q$-binomial coefficient is defined as follows 
(see, for example, \cite{GasparRahman}):
\be
\label{qbin}
{\,n\,\atopwithdelims []\,k\,}_{q}
:=\frac{[n]_{q}!}{[k]_{q}![n-k]_{q}! }
=\frac{(q^{-n};q)_k}{(q;q)_k}(-1)^kq^{nk-({}^k_2)} \qquad 
(0\leqq k\leqq n).
\ee 

The basic (or $q$-) hypergeometric function 
of the variable $z$ and with $\mathfrak{r}$ numerator 
and $\mathfrak{s}$ denominator parameters 
is defined as follows (see, for details, the monographs by 
Slater \cite[Chapter 3]{Slater} 
and by Srivastava and Karlsson 
\cite[p. 347, Eq. (272)]{SrivastavaKarlsson}; 
see also \cite{HMS-IMAJAM1983-1984} and \cite{Koekock}): 

$${}_{\mathfrak r}\Phi_{\mathfrak s}\left[
\begin{array}{rr}
a_1, a_2,\cdots, a_{\mathfrak r};\\
\\
b_1,b_2,\cdots,b_{\mathfrak s};
\end{array}\,
q;z\right]
:=\sum_{n=0}^\infty\Big[(-1)^n \; 
q^{\binom{n}{2}}\Big]^{1+{\mathfrak s}-{\mathfrak r }}
\,\frac{(a_1, a_2,\cdots, a_{\mathfrak r};q)_n}
{(b_1,b_2,\cdots,b_{\mathfrak s};q)_n}
\; \frac{z^n}{(q;q)_n},$$
where $q\neq 0$ when ${\mathfrak r }>{\mathfrak s}+1$. 
We also note that

$${}_{\mathfrak r+1}\Phi_{\mathfrak r}\left[
\begin{array}{rr}
a_1, a_2,\cdots, a_{\mathfrak r+1}\\
\\
b_1,b_2,\cdots,b_{\mathfrak r };
\end{array}\,
q;z\right]
=\sum_{n=0}^\infty \frac{(a_1, a_2,\cdots, a_{\mathfrak r+1};q)_n}
{(b_1,b_2,\cdots,b_{\mathfrak r};q)_n}\;\frac{ z^n}{(q;q)_n}.$$

Here, in our present investigation, we are mainly concerned  
with the Cauchy polynomials $p_n(x,y)$ as given
below (see \cite{Chen2003} and \cite{GasparRahman}):
\bea
\label{def}
p_n(x,y):=(x-y)(x- qy)\cdots (x-q^{n-1}y)
=\left(\frac{y}{x};q\right)_n\,x^n 
\eea
together with the following Srivastava-Agarwal 
type generating function 
(see also \cite{Cao-Srivastava2013}):  
\be
\label{Srivas}
\sum_{n=0}^\infty  
p_n (x,y)\frac{(\lambda;q)_n\,t^n}{(q;q)_n}
={}_{2}\Phi_1\left[
\begin{array}{rr}
\lambda,y/x;\\
\\
0;
\end{array}\,
q;  xt\right].
\ee

For $\lambda=0$ in (\ref{Srivas}), we get the 
following simpler generating function \cite{Chen2003}:
\be
\label{gener}
\sum_{n=0}^{\infty} p_n(x,y)
\frac{t^n }{(q;q)_n} = 
\frac{(yt;q)_\infty}{(xt;q)_\infty}.
\ee

The generating function (\ref{gener}) 
is also the homogeneous version
of the Cauchy identity or the following 
$q$-binomial theorem (see, for example, \cite{GasparRahman}, 
\cite{Slater} and \cite{SrivastavaKarlsson}):
\be
\label{putt}
\sum_{k=0}^{\infty} 
\frac{(a;q)_k }{(q;q)_k}\;z^{k}={}_{1}\Phi_0\left[
\begin{array}{rr}
a;\\
\\
\overline{\hspace{3mm}}\,;
\end{array} \,
q;z\right]=\frac{(az;q)_\infty}{(z;q)_\infty}\qquad (|z|<1). 
\ee
Upon further setting $a=0$, this last relation (\ref{putt}) 
becomes Euler's identity  
(see, for example, \cite{GasparRahman}):
\be
\label{q-expo-alpha}
\sum_{k=0}^{\infty} \frac{z^{k}}{(q;q)_k}=\frac{1}{(z;q)_\infty}
\qquad (|z|<1)
\ee
and its inverse relation given below \cite{GasparRahman}:
\be
\label{q-Expo-alpha}
 \sum_{k=0}^{\infty}\frac{(-1)^k
}{(q;q)_k}\; q^{\binom{k}{2}}\,z^{k}=(z;q)_\infty.
\ee

The Jackson's $q$-difference or $q$-derivative operators $D_q$ and $D_{q^{-1}}$ 
are defined as follows (see, for example, \cite{Jackson1908,VKac2002,GasparRahman}): 
\be 
\label{dd2}
 {D}_q\big\{f(x)\big\}:=\frac{f(x)-f( q x)}{ (1-q) x} \quad \text{and} \quad 
D_{q^{-1}}\big\{f(x)\}:=\frac{f(q^{-1}x)-f(x)}{(q^{-1}-1)x}.
\ee
Evidently, in the limit when $q\rightarrow 1-$, we have
$$\lim_{q\rightarrow 1-} 
\left\{ {D}_q\big\{f(x)\big\} \right\}
=f'(x)\quad \text{and} \quad \lim_{q\rightarrow 1-} 
\left\{ D_{q^{-1}}\big\{f(x)\big\} \right\}
=f'(x),$$
provided that the derivative $f'(x)$ exists.   \\

Suppose that $D_q$ acts on the variable $a$. Then
we have the $q$-identities asserted by Lemma \ref{MALM} below.

\begin{lemma}
\label{MALM}
Each of the following $q$-identities holds true for
the $q$-derivative operator $D_q$ acting on the variable $a$$:$
\begin{align}
\label{id1}
{D}_q^k \left\{\frac{1}{(as;q)_\infty}\right\}
=\frac{\big[(1-q)^{-1}s\big]^k}{(as;q)_\infty},
\end{align}
\begin{align}
\label{id2}
\big(D_{q^{-1}}\big)^k \left\{\frac{1}{(as;q)_\infty}\right\}
=\frac{q^{-\binom{k}{2}}}{(asq^{-k};q)_\infty}\;
\; \left(\frac{s}{1-q}\right)^k,
\end{align}
\begin{align}
\label{id3}
{D}_q^k \left\{{(as;q)_\infty}\right\}
=(-1)^k \;q^{\binom{k}{2}}\; 
(asq^k;q)_\infty\; \left(\frac{s}{1-q}\right)^k
\end{align}
and
\begin{align} 
\label{id4}
\big(D_{q^{-1}}\big)^k \left\{{(as;q)_\infty}\right\}
=(as;q)_\infty\;
\left(-\frac{s}{1-q}\right)^k.
\end{align}
\end{lemma}

The Leibniz rules for the $q$-derivative operators 
$D_q$ and $D_{q^{-1}}$ are given by 
(see, for example, \cite{Liu97} and \cite{Liu98})
\be
\label{Lieb}
D_q^n\left\{f(x)g(x)\right\}=\sum_{k=0}^n q^{k(k-n)}\;
\qbinomial{n}{k}{q} \;
D_q^k\left\{f(x)\right\} \; D_q^{n-k}\left\{g(q^kx)\right\}
\ee 
and
\be
\label{Lieb1}
D_{q^{-1}}^n\left\{f(x)g(x)\right\}=\sum_{k=0}^n 
\qbinomial{n}{k}{q} \; 
D_{q^{-1}}^k\left\{f(x)\right\} \;
D_{q^{-1}}^{n-k}\left\{g(q^{-k}x)\right\},
\ee
where $D_q^0$ and $D_{q^{-1}}^0$ are 
understood to be the identity operators.

\begin{lemma} Suppose that $q$-difference operator 
$D_q$ acts on the variable $a$. Then
\label{dAMM}
\begin{equation} 
D_q^n\left\{\frac{(as;q)_\infty}{(a\omega;q)_\infty}\right\}
=\left(\frac{\omega}{1-q}\right)^n\; \frac{(s/\omega;q)_n}
{(as;q)_n}\;\frac{(as;q)_\infty}{( a\omega;q)_\infty}\label{aberll}
\end{equation}
and
\begin{equation}
D_{q^{-1}}^n\left\{\frac{(as;q)_\infty}{(a\omega;q)_\infty}\right\}  
=\left(-\frac{q}{(1-q)a}\right)^n\; \frac{(s/\omega;q)_n}{(q/(a\omega);q)_n}   
\; \frac{(as;q)_\infty}{(a\omega;q)_\infty}.   \label{abell}
\end{equation} 
\end{lemma}

\begin{proof}
Suppose that the operator 
${D}_{q}$  acts upon the variable $a$.
In light of (\ref{Lieb}), we then find that
\begin{align}\label{proof1}
D_q^n\left\{\frac{(as;q)_\infty}{(a\omega;q)_\infty}\right\}
&=\sum_{k=0}^n q^{k(k-n)}\qbinomial{n}{k}{q}\;
D_q^k\left\{(as;q)_\infty \right\}\; 
D_q^{n-k}\left\{\frac{1}{(a\omega q^{k};q)_\infty}\right\}\notag\\     
&=\sum_{k=0}^n q^{k(k-n)+\binom{k}{2}}\;\qbinomial{n}{k}{q}\;
(asq^k;q)_\infty\;\left(-\frac{s}{1-q}\right)^k\;
\frac{\big[(1-q)^{-1}\omega q^{k}\big]^{n-k}}{(a\omega q^{k};q)_\infty}\notag \\    
&=\frac{(as;q)_\infty}{(a\omega;q)_\infty}\;\left(\frac{\omega}{1-q}\right)^n\;
\sum_{k=0}^n \frac{(q^{-n},a\omega;q)_k}{(as,q;q)_k}\;
\left(\frac{s q^n}{\omega}\right)^k \notag \\
&=\left(\frac{\omega}{1-q}\right)^n\;\frac{(as;q)_\infty}{(a\omega;q)_\infty}\;
{}_{2}\Phi_{1}\left[
\begin{array}{rr}q^{-n}, a\omega;\\
\\ 
as; 
\end{array} \,
q; \frac{s q^n}{\omega}\right],
\end{align}
where we have appropriately applied the formulas (\ref{id1}), 
(\ref{id3}) and (\ref{qbin}).
The proof of the first assertion (\ref{aberll}) of Lemma \ref{dAMM} 
is completed by using the relation (\ref{chuv}) in (\ref{proof1}). \\

Similarly, by using the relation (\ref{Lieb1}), we have
\begin{align}\label{proof2}
D_{q^{-1}}^n\left\{\frac{(as;q)_\infty}{(a\omega;q)_\infty}
\right\}&=\sum_{k=0}^n \qbinomial{n}{k}{q}\;  D_{q^{-1}}^k
\left\{\frac{1}{(a\omega ;q)_\infty}\right\}\; D_{q^{-1}}^{n-k}
\left\{(asq^{-k};q)_\infty \right\}\notag \\
&=\sum_{k=0}^n \qbinomial{n}{k}{q}\; \frac{q^{-\binom{k}{2}}
\;\left(\frac{\omega}{1-q}\right)^k}
{(a\omega q^{-k};q)_\infty}
(-1)^{n-k}\;\left(\frac{sq^{-k}}{1-q}\right)^{n-k}\;(asq^{-k};q)_\infty
\notag \\
&=\frac{(as;q)_\infty}{(a\omega;q)_\infty}\left(-\frac{s}{1-q}\right)^n
\sum_{k=0}^n \qbinomial{n}{k}{q} \;
\frac{(-1)^k\; q^{-\binom{k}{2}+k(k-n)}\;(as q^{-k};q)_k}
{(a\omega q^{-k};q)_k}\left(\frac{\omega}{s}\right)^k \notag \\
&=\frac{(as;q)_\infty}{(a\omega;q)_\infty}\; \left(-\frac{s}{1-q}\right)^n\; \sum_{k=0}^n
\;\frac{\big(q^{-n},q/(as);q\big)_k \,q^k}{\big(q,q/(a\omega);q\big)_k}\notag \\
&=\left(-\frac{s}{1-q}\right)^n \;\frac{(as;q)_\infty}{(a\omega;q)_\infty}\;{}_{2}\Phi_{1}
\left[
\begin{array}{rr}
q^{-n},q/(as);\\
\\
q/(a\omega);
\end{array} 
\;q;q\right],
\end{align}
where we have appropriately used the relation (\ref{qbin}).

Finally, by using the relation (\ref{qchuv}) in (\ref{proof2}), 
we are led to the second assertion (\ref{abell}) 
of Lemma \ref{dAMM}. 
We thus have completed the proof of Lemma \ref{dAMM}. 
\end{proof}

\begin{remark} \label{rem1.1}
{\rm For $s=0$ and $\omega=s$, the assertions $(\ref{aberll})$ and 
$(\ref{abell})$ of Lemma $\ref{dAMM}$ reduce to the assertions $(\ref{id1})$ 
and $(\ref{id2})$ of Lemma $\ref{MALM}$.
Moreover, for $\omega=0$, the assertions $(\ref{aberll})$ and $(\ref{abell})$ 
of Lemma $\ref{dAMM}$ reduce to the assertions $(\ref{id3})$ and $(\ref{id4})$  
of Lemma $\ref{MALM}$.}
\end{remark}

\begin{lemma}{\rm (see, for example, 
\cite[Eq. (0.58) and Eq. (II.6)]{GasparRahman})}
The $q$-Chu-Vandermonde formulas are given by
\be
\label{chuv}
{}_{2}\Phi_1\left[
\begin{array}{rr} q^{-n}, a;\\
\\ 
c;
\end{array}\,
q; \frac{cq^n}{a} \right]=\frac{(c/a;q)_n}{(c;q)_n} 
\qquad (n\in \mathbb{N}_0:=\mathbb{N}\cup\{0\})
\ee
and
\be
\label{qchuv}
{}_{2}\Phi_1\left[
\begin{array}{rr} q^{-n}, a;\\
\\
c;
\end{array}\,
q; q \right]=\frac{(c/a;q)_n}{(c;q)_n}\;a^n
\qquad (n\in \mathbb{N}_0:=\mathbb{N}\cup\{0\}).
\ee
\end{lemma}
 
We now state and prove the $q$-difference formulas asserted by
Theorem \ref{Zition} below. 

\begin{theorem}  
\label{Zition} Suppose that the $q$-difference operators 
${D}_{q}$ and $D_{q^{-1}}$ act upon the variable $a$. Then
\begin{equation}
D_q^n\left\{\frac{(as;q)_\infty}{(at,a\omega;q)_\infty}\right\} 
=\left(\frac{t}{1-q}\right)^n \;\frac{(as;q)_\infty}{(at,a\omega;q)_\infty}  
\; {}_{3}\Phi_{1}\left[
\begin{array}{rr}
q^{-n}, s/\omega, at;\\
\\
as;
\end{array} 
\;q;\frac{\omega q^n}{t}\right],\label{oker1}
\end{equation}
and

\begin{equation}
D_{q^{-1}}^n\left\{\frac{(as,at;q)_\infty}{(a\omega;q)_\infty}\right\} 
=\left(-\frac{t}{1-q}\right)^n\; \frac{(at,as;q)_\infty}{(a\omega;q)_\infty} 
\; {}_{3}\Phi_{2}\left[
\begin{array}{rr}
q^{-n}, q/(at),s/\omega;\\ 
\\
q/(a\omega),0;
\end{array} 
\,q; q\right].\label{oker}
\end{equation}
\end{theorem}

\begin{proof}
Suppose first that the $q$-difference operator 
${D}_{q}$ acts upon the variable $a$.
Then, in light of (\ref{Lieb}), and by using the relations 
(\ref{aberll}) and (\ref{id1}), it is easily seen that
\begin{align}\label{proof3}
D_q^n\left\{\frac{(as;q)_\infty}{(a\omega,at;q)_\infty}\right\}
&=\sum_{k=0}^n \qbinomial{n}{k}{q}\; q^{k(k-n)}\; D_q^k\left\{\frac{(as;q)_\infty}
{(a\omega;q)_\infty}\right\} D_q^{n-k}
\left\{\frac{1}{(atq^{k};q)_\infty}\right\}\notag \\
&=\sum_{k=0}^n \qbinomial{n}{k}{q}\;q^{k(k-n)}
\; \frac{(s/\omega;q)_k}{(as;q)_k}\;\frac{(as;q)_\infty}{(a\omega;q)_\infty}
\; \frac{[(1-q)^{-1}tq^{k}]^{n-k}}{(atq^{k};q)_\infty}\; 
\left(\frac{\omega}{1-q}\right)^k\notag \\  
&= \frac{(as;q)_\infty}{(a\omega,at;q)_\infty}\;
\left(\frac{t}{1-q}\right)^n\; \sum_{k=0}^n\; 
\frac{(-1)^k\;q^{-\binom{k}{2}}\; 
(q^{-n},s/\omega,at;q)_k}{(as,q;q)_k}
\left(\frac{\omega q^n}{t}\right)^k
\notag \\
&=\left(\frac{t}{1-q}\right)^n\;
\frac{(as;q)_\infty}{(a\omega,at;q)_\infty}  
\; {}_{3}\Phi_{1}\left[
\begin{array}{rr}
q^{-n}, s/\omega, at;\\
\\
as;
\end{array} \,
q; \frac{\omega q^n}{t}\right].
\end{align}

Similarly, by using (\ref{Lieb1}), 
we find for the $q$-difference operator 
$D_{q^{-1}}$ acting on the variable $a$ that
\begin{align} \label{proof4}
D_{q^{-1}}^n\left\{\frac{(at,as;q)_\infty}
{(a\omega;q)_\infty}\right\}
&=\sum_{k=0}^n \qbinomial{n}{k}{q} \; D_{q^{-1}}^k
\left\{\frac{(as;q)_\infty}{(a\omega;q)_\infty}\right\} \;
D_{q^{-1}}^{n-k}\left\{(atq^{-k};q)_\infty \right\}\notag \\
&=\sum_{k=0}^n \qbinomial{n}{k}{q}\;
\frac{\big[-qa^{-1}(1-q)^{-1}\big]^k\;
(s/\omega;q)_k(as;q)_\infty}
{\big(q/(a\omega);q\big)_k(a\omega;q)_\infty}
\;  (atq^{-k};q)_\infty \; 
\left(-\frac{tq^{-k}}{1-q}\right)^{n-k}
\notag \\
&=\left(-\frac{t}{1-q}\right)^n\;\frac{(as,at;q)_\infty}
{(a\omega;q)_\infty}\sum_{k=0}^n \qbinomial{n}{k}{q}\;
q^{k(1+k-n)} \;(ta)^{-k}\; \frac{(s/\omega,atq^{-k};q)_k}
{\big(q/(a\omega);q\big)_k} \notag \\
&=\left(-\frac{t}{1-q}\right)^n\; \frac{(as,at;q)_\infty}
{(a\omega;q)_\infty}\;\sum_{k=0}^n
\frac{\big(q^{-n},s/\omega,q/(at);q\big)_k}
{\big(q,q/(a\omega);q\big)_k}\;q^{k },
\end{align}
where we have appropriately used the relation (\ref{qbin}) as well.

Equations (\ref{proof3}) and (\ref{proof4}), together, 
complete the proof of Theorem \ref{Zition}.
\end{proof}

\begin{remark} \label{rem1.2}
{\rm Upon first setting $\omega=0,$ we put $t=\omega$ 
in the assertion $(\ref{oker1})$ of Theorem $\ref{Zition},$  
Then, if we make use of the identity $(\ref{chuv}),$
we get $(\ref{aberll})$. Furthermore, upon setting $s=0$ in the 
assertion $(\ref{oker})$ of Theorem $\ref{Zition},$ if we make use of 
the $q$-Chu-Vandermonde  formula $(\ref{qchuv}),$ 
we get $(\ref{abell})$.}
\end{remark}
 
This paper is organized as follows. 
In Section \ref{qdifff}, we introduce two homogeneous $q$-difference   
operators $\mathbb{T}({\bf a},{\bf b},cD_{q})$ and 
$\mathbb{E}({\bf a},{\bf b},cD_{q^{-1}})$. In addition, we define two 
families of the generalized  Al-Salam-Carlitz $q$-polynomials 
$\phi_n^{({\bf a},{\bf b})}(x,y|q)$ and $\psi_n^{({\bf a},{\bf b})}(x,y|q)$ 
and represent each of the families in terms of the homogeneous 
$q$-difference operators $\mathbb{T}({\bf a},{\bf b},cD_{q})$ and 
$\mathbb{E}({\bf a},{\bf b},cD_{q^{-1}})$. We also derive generating functions
for these families of the generalized  Al-Salam-Carlitz $q$-polynomials. 
In Section \ref{section3}, we first give the Rogers type
formulas and the extended Rogers type formulas.  
The Srivastava-Agarwal type generating functions involving the generalized  
Al-Salam-Carlitz $q$-polynomials are derived in Section \ref{qdiffff}.
Finally, in our last section (Section \ref{Conclusion}), we present the
concluding remarks and observations concerning our present investigation.

\section{Generalized Al-Salam-Carlitz $q$-Polynomials}
\label{qdifff}

In this section, we first introduce 
two homogeneous $q$-difference operators
$\mathbb{T}({\bf a},{\bf b},cD_{q})$ and 
$\mathbb{E}({\bf a},{\bf b},cD_{q^{-1}})$
which are defined by

\begin{equation}
\label{operator0}
\mathbb{T}({\bf a},{\bf b}, cD_{q}):=\sum_{k=0}^\infty\;
\frac{(a_1,a_2,\cdots, a_{\mathfrak r+1};q)_k}
{(q,b_1,b_2,\cdots,b_{\mathfrak r};q)_k}(cD_{q})^k 
\end{equation}
and
\begin{equation}
\label{operator1}
\mathbb{E}({\bf a},{\bf b},cD_{q^{-1}}):=\sum_{k=0}^\infty\;
\frac{(a_1,a_2,\cdots,a_{\mathfrak r+1};q)_k}
{(q,b_1,b_2,\cdots,b_{\mathfrak r};q)_k}(cD_{q^{-1}})^k, 
\end{equation}
where, for convenience, 
$${\bf a}=(a_1,a_2,\cdots,a_{\mathfrak r+1})
\qquad \text{and} \qquad 
{\bf b}=(b_1,b_2,\cdots,b_{\mathfrak r}).$$

We now derive the identities (\ref{MaA}) and (\ref{MaA1}) below,  
which will be used later in order to derive the generating functions, 
the Rogers type formulas, the extended Rogers type formulas and 
the Srivastava-Agarwal type generating functions involving   
the families of the generalized Al-Salam-Carlitz $q$-polynomials.

\begin{theorem}
\label{prei}  
Suppose that the $q$-difference operator $D_q$ 
acts on the variable $a$. Then
\begin{align}
\label{MaA}
\mathbb{T}({\bf a},{\bf b},cD_{q})\left\{\frac{(as;q)_\infty} 
{(a\omega,at;q)_\infty}\right\}
&=\frac{(as;q)_\infty}{(a\omega,at;q)_\infty}\;\sum_{n=0}^\infty 
\;\frac{(a_1,a_2,\cdots,a_{\mathfrak r+1};q)_n}
{(q,b_1,b_2,\cdots,b_{\mathfrak r};q)_n}\; 
\left(\frac{ct}{1-q}\right)^n\notag \\
&\qquad \cdot \, {}_{3}\Phi_{1}\left[
\begin{array}{rr}
q^{-n},s/\omega, at;\\ 
\\
as;
\end{array}\,
q; \frac{\omega q^n}{t}\right] \qquad (\max\{|a\omega|,|at|\} <1)
\end{align}
and

\begin{align}
\label{MaA1}
\mathbb{E}({\bf a},{\bf b},cD_{q^{-1}})
\left\{\frac{(at,as;q)_\infty}{(a\omega;q)_\infty}\right\}
&=\frac{(at,as;q)_\infty}{(a\omega;q)_\infty}\;
\sum_{n=0}^\infty \frac{(a_1,a_2,\cdots,a_{\mathfrak r+1};q)_n}
{(q,b_1,b_2,\cdots,b_{\mathfrak r };q)_n}\; 
\left(-\frac{ct}{1-q}\right)^n\notag \\
&\qquad \cdot\, {}_{3}\Phi_{2}\left[
\begin{array}{rr}
q^{-n}, q/(at),s/\omega;\\ 
\\
q/(a\omega),0;
 \end{array} \,
q; q\right] \qquad (|a\omega|<1).
\end{align}
\end{theorem}

\begin{proof}
Suppose that the operators $D_q$ and $D_{q^{-1}}$ act on the variable $a$.
We observe by applying (\ref{oker1}) that
\begin{align}\label{proof5}
&\mathbb{T}({\bf a},{\bf b},cD_{q})\left\{\frac{(as;q)_\infty} 
{(a\omega,at;q)_\infty}\right\}\notag \\
&\qquad =\sum_{n=0}^\infty\;\frac{(a_1,a_2,\cdots,a_{\mathfrak r+1};q)_n\,c^n}
{(q,b_1,b_2,\cdots,b_{\mathfrak r };q)_n}\;D_{q}^n\left\{\frac{(as;q)_\infty} 
{(a\omega,at;q)_\infty}\right\}\notag \\
&\qquad =\sum_{n=0}^\infty\;\frac{ (a_1,a_2,\cdots,
a_{\mathfrak r+1};q)_n\,c^n}
{(q,b_1,b_2,\cdots, b_{\mathfrak r};q)_n}\;\frac{(as;q)_\infty}
{(a\omega,at;q)_\infty}\; 
\left(\frac{t}{1-q}\right)^n \;{}_{3}\Phi_{1}\left[
\begin{array}{rr}
q^{-n}, s/\omega, at;\\
\\
as;
\end{array} \,
q; \frac{\omega q^n}{t}\right]
\notag\\
&\qquad=\frac{(as;q)_\infty}{(a\omega,at;q)_\infty}\;
\sum_{n=0}^\infty\;
\frac{(a_1,a_2,\cdots,a_{\mathfrak r+1};q)_n} 
{(q,b_1,b_2,\cdots, b_{\mathfrak r};q)_n}\;
{}_{3}\Phi_{1}\left[
\begin{array}{rr}
q^{-n}, s/\omega, at;\\
\\
as;
\end{array}\, 
q; \frac{\omega q^n}{t}\right]\; 
\left(\frac{ct}{1-q}\right)^n. 
\end{align}

Similarly, by applying (\ref{oker}), we find that 
\begin{align}\label{proof6}
&\mathbb{E}({\bf a},{\bf b},cD_{q^{-1}})\left\{\frac{(at,as;q)_\infty} 
{(a\omega;q)_\infty}\right\}\notag \\
&\qquad =\sum_{n=0}^\infty\frac{(a_1,a_2,\cdots, a_{\mathfrak r+1};q)_n\,(c)^n}
{(q,b_1,b_2,\cdots,b_{\mathfrak r};q)_n}\;
D_{q^{-1}}^n\left\{\frac{(at,as;q)_\infty}{(a\omega;q)_\infty}\right\} \notag\\ 
&\qquad=\sum_{n=0}^\infty \;\frac{(a_1,a_2,\cdots,a_{\mathfrak r+1};q)_n \,
(c)^n}{(q,b_1,b_2,\cdots,b_{\mathfrak r};q)_n}\;
\frac{(at,as;q)_\infty}{(a\omega;q)_\infty}\left(\frac{t}{ q-1}\right)^n \;
{}_{3}\Phi_{2}\left[
\begin{array}{rr}
q^{-n}, q/(at),s/\omega;\\ 
\\
q/(a\omega),0;
 \end{array}\, 
q; q\right]
\notag \\
&\qquad=\frac{(at,as;q)_\infty}{(a\omega;q)_\infty}\;
\sum_{n=0}^\infty \frac{(a_1,a_2,\cdots, a_{\mathfrak r+1};q)_n}
{(q,b_1,b_2,\cdots, b_{\mathfrak r};q)_n}\; 
\left(\frac{ct}{q-1}\right)^n\; 
{}_{3}\Phi_{2}\left[
\begin{array}{rr}
q^{-n}, q/(at),s/\omega;\\ 
\\
q/(a\omega),0;
\end{array} \,
q; q\right], 
\end{align}
as asserted by Theorem \ref{prei}.
\end{proof}

\noindent 
\textbf{Definition.}
{\rm In terms of the $q$-binomial coefficient, the families of the  
generalized Al-Salam-Carlitz $q$-polynomials 
$\phi_n^{({\bf a},{\bf b})}(x,y|q)$ and 
$\psi_n^{({\bf a},{\bf b})}(x,y|q)$
are defined by
\begin{equation}
\label{Carl}
\phi_n^{({\bf a}, {\bf b})}(x,y|q)
=\sum_{k=0}^n\qbinomial{n}{k}{q}\;
\frac{(a_1,a_2,\cdots,a_{\mathfrak r+1};q)_k}
{(b_1,b_2,\cdots,b_{\mathfrak r};q)_k}\; x^k\;y^{n-k}
\end{equation}
and
\begin{equation}
\label{Carl1}
\psi_n^{({\bf a},{\bf b})}(x,y|q)
=\sum_{k=0}^n \qbinomial{n}{k}{q}\;
\frac{(a_1,a_2,\cdots,a_{\mathfrak r+1};q)_k}
{(b_1,b_2,\cdots,b_{\mathfrak r};q)_k} 
q^{\binom{k+1}{2}-nk}\,x^k\;y^{n-k}.
\end{equation}
}

\medskip

\noindent
\textbf{Proposition.}
{\it Suppose that the operators $D_q$ and $D_{q^{-1}}$ 
act on the variable $y$. Then
\begin{equation}
\label{Propo12}
\phi_n^{({\bf a},{\bf b})}(x,y|q)
=\mathbb{T}\big({\bf a},{\bf b},(1-q)xD_{q}\big)\{y^n\}
\quad {\it and} \quad \psi_n^{({\bf a},{\bf b})}(x,y|q)
=\mathbb{E}\big({\bf a},{\bf b},-(1-q)xD_{q^{-1}}\big)\{y^n\}
\end{equation}
in terms of the operators $(\ref{operator0})$ and $(\ref{operator1})$}.

\begin{theorem}
{\rm $\big[$Generating function for $\phi_n^{({\bf a},
{\bf b})}(x,y,z|q)$ and $\psi_n^{({\bf a},{\bf b})}(x,y,z|q)\big]$}
\label{T1}
Each of the following generating functions holds true$:$
\bea
\label{gen}
\sum_{n=0}^\infty \phi_n^{({\bf a},{\bf b})}(x,y,z|q)\; \frac{t^n}{(q;q)_n}
=\frac{1}{(yt;q)_\infty}\;{}_{{\mathfrak r+1}}\Phi_{{\mathfrak r}}\left[
\begin{array}{rr}
a_1,a_2,\cdots,a_{\mathfrak r+1};\\ 
\\
b_1,b_2,\cdots,b_{\mathfrak r}; 
\end{array}\, 
q; xt\right] \qquad  (\max\{|xt|,|yt|\}<1)
\eea
and
\bea
\label{gen1}
\sum_{n=0}^\infty (-1)^n\; q^{\binom{n}{2}} \;
\psi_n^{({\bf a},{\bf b})}(x,y|q)\; \frac{ 
t^n}{(q;q)_n}=(yt;q)_\infty \;{}_{{\mathfrak r+1}}\Phi_{{\mathfrak r}}\left[
\begin{array}{rr}
a_1,a_2,\cdots,a_{\mathfrak r+1};\\
\\
b_1,b_2,\cdots,b_{\mathfrak r};
 \end{array}\,
q; xt\right] \qquad (|xt|<1).
\eea
\end{theorem}

In our proof of Theorem \ref{T1}, the following 
easily derivable Lemma will be needed.

\begin{lemma} \label{lem2.1}
Suppose that the operators $D_q$ and $D_{q^{-1}}$ 
act on the variable $a$. Then
\be 
\label{Ma2}
\mathbb{T}({\bf a},{\bf b},cD_{q})
\left\{\frac{1}{(as;q)_\infty}\right\} 
=\frac{1}{(as;q)_\infty}\; 
{}_{{\mathfrak r+1}}\Phi_{{\mathfrak r+1}}
\left[
\begin{array}{rr}
a_1,a_2,\cdots,a_{\mathfrak r+1};\\ 
\\
b_1,b_2,\cdots,b_{\mathfrak r}; 
\end{array} \,
q;\frac{cs}{1-q}\right]
\ee
$$\left(\max\left\{|as|,\left|\frac{cs}{1-q}\right|\right\}<1\right).$$
and

\be
\label{Ma5}
\mathbb{E}({\bf a},{\bf b},-cD_{q^{-1}})
\left\{(as;q)_\infty \right\}
=(as;q)_\infty \;{}_{{\mathfrak r+1}}
\Phi_{{\mathfrak r}}\left[
\begin{array}{rr}
a_1,a_2,\cdots, a_{\mathfrak r+1};\\ 
\\
b_1,b_2,\cdots, b_{\mathfrak r}; 
\end{array} \,
q;\frac{cs}{1-q}\right] \qquad 
\left(\left|\frac{cs}{1-q}\right|<1\right).
\ee
\end{lemma}

\begin{proof}[Proof of Theorem $\ref{T1}$] 
We suppose that the $q$-difference operator 
${D}_{q}$ acts upon the variable $y$. 
In light of the formulas in (\ref{Propo12}),  
and by applying (\ref{Ma2}), it is readily seen that 
\begin{align}\label{proof7}
\sum_{n=0}^\infty \phi_n^{({\bf a},{\bf b})}(x,y|q)
\;\frac{t^n}{(q;q)_n}
&=\sum_{n=0}^\infty 
\mathbb{T}\big({\bf a},{\bf b},(1-q)xD_q\big)
\left\{y^n\right\}\; \frac{t^n}{(q;q)_n}\notag \\
&=\mathbb{T}\big({\bf a},{\bf b},(1-q)xD_q\big)
\left\{\sum_{n=0}^\infty \;
\frac{(yt)^n}{(q;q)_n}\right\}\notag \\
&=\mathbb{T}\big({\bf a},{\bf b},(1-q)xD_q\big)
\left\{\frac{1}{(yt;q)_\infty}\right\}\notag \\
&=\frac{1}{(yt;q)_\infty}\;{}_{{\mathfrak r+1}}
\Phi_{{\mathfrak r }}\left[
\begin{array}{rr}
a_1,a_2,\cdots, a_{\mathfrak r+1};\\
\\
b_1,b_2,\cdots, b_{\mathfrak r};
\end{array} \,
q; xt\right].
\end{align}

Similarly, we have 
\begin{align}\label{proof8}
\sum_{n=0}^\infty (-1)^n\;q^{\binom{n}{2}}\;
\psi_n^{({\bf a},{\bf b})}(x,y|q)\;\frac{ 
t^n}{(q;q)_n}&= \sum_{n=0}^\infty 
\mathbb{E}\big({\bf a},{\bf b},-(1-q)xD_{q^{-1}}\big)
\left\{y^n\right\}\;(-1)^n\;q^{\binom{n}{2}}\; 
\frac{t^n}{(q;q)_n}\notag \\
&=\mathbb{E}\big({\bf a},{\bf b},-(1-q)xD_{q^{-1}}\big)
\left\{\sum_{n=0}^\infty\;
\frac{(-1)^n\; q^{\binom{n}{2}}(yt)^n}{(q;q)_n}\right\}\notag \\
&=\mathbb{E}\big({\bf a},{\bf b},-(1-q)xD_{q^{-1}}\big)
\left\{(yt;q)_\infty \right\}.
\end{align}
The proof of Theorem \ref{T1} can now be completed by making use 
of the relation (\ref{Ma5}).
\end{proof}

\section{The Rogers Type Formulas and the Extended Rogers Type Formulas}
\label{section3}
In this section, we use the assertions in (\ref{Propo12}) to derive several 
$q$-identities such as the Rogers type formulas and the extended Rogers type formulas  
for the families  of the generalized Al-Salam-Carlitz $q$-polynomials  
$\phi_{n}^{({\bf a},{\bf b})}(x,y|q)$ and 
$\psi_{n}^{({\bf a},{\bf b})}(x,y|q)$.

\begin{theorem}
{\rm $\big[$Rogers type formula for 
$\phi_{n}^{({\bf a},{\bf b})}(x,y|q)\big]$}
\label{exgend}
The following Rogers type formula
holds true for $\phi_{n}^{({\bf a},{\bf b})}(x,y|q)$$:$ 
\begin{align}
\label{exgen}
&\sum_{n=0}^\infty\sum_{m=0}^\infty   
\phi_{n+m}^{({\bf a},{\bf b})}(x,y|q)\; \frac{t^n}{(q;q)_n}\;
\frac{s^m}{(q;q)_m} \notag \\ 
&\qquad =\frac{1}{(yt,ys;q)_\infty}\;\sum_{n=0}^\infty 
\;\frac{(a_1,a_2,\cdots,a_{\mathfrak r+1};q)_n}
{(q,b_1,b_2,\cdots, b_{\mathfrak r };q)_n}\; (xt)^n\; 
{}_{2}\Phi_{0}\left[
\begin{array}{rr}
q^{-n},yt;\\
 \\
\overline{\hspace{9mm}}\;;
\end{array}\,
q;\frac{s q^n}{t}\right] 
\end{align}
$$(\max\{|yt|,|ys|\}<1).$$
\end{theorem}
\begin{theorem}
{\rm $\big[$Rogers type formula for 
$\psi_{n}^{({\bf a}, {\bf b})}(x,y|q)\big]$}
\label{exgend1}
The following Rogers type formula
holds true for $\psi_{n}^{({\bf a},{\bf b})}(x,y|q)$$:$ 
\begin{align}
\label{exgen1}
&\sum_{n=0}^\infty\sum_{m=0}^\infty\; 
(-1)^{n+m}\;q^{\binom{n}{2}+\binom{m}{2}}\;    
\psi_{n+m}^{({\bf a},{\bf b})}(x,y|q) \;
\frac{t^n }{(q;q)_n}\; \frac{s^m}{(q;q)_m}\notag \\
&\qquad =(yt,ys;q)_\infty
\;\sum_{n=0}^\infty\; 
\frac{(a_1,a_2,\cdots,a_{\mathfrak r+1};q)_n}
{(b_1,b_2,\cdots,b_{\mathfrak r};q)_n}\;(xt)^n\; 
{}_{2}\Phi_{1}\left[
\begin{array}{rr}
q^{-n}, q/(yt);\\  
\\
0;
\end{array}\,
q; ys\right].
\end{align}
\end{theorem}

In order to prove Theorems \ref{exgend} and \ref{exgend1}, 
we need Lemma \ref{lem3.1} below. 

\begin{lemma}\label{lem3.1} 
It is asserted that
\begin{align}
\label{Ma2A}
&\mathbb{T}({\bf a},{\bf b},cD_{q})
\left\{\frac{1}{(a\omega,at;q)_\infty}\right\}\notag \\
&\qquad =\frac{(as;q)_\infty}{(a\omega,at;q)_\infty}\;
\sum_{n=0}^\infty \;\frac{(a_1,a_2,\cdots,a_{\mathfrak r+1};q)_n}
{(q,b_1,b_2,\cdots, b_{\mathfrak r };q)_n}
\;       
{}_{2}\Phi_{0}\left[
\begin{array}{rr}
q^{-n},at;\\
\\
\overline{\hspace{9mm}}\;;
\end{array}\,
q; \frac{\omega q^n}{t}\right] 
\; \left(\frac{ct}{1-q}\right)^n
\end{align}
$$(\max\{|a\omega|,|at|\}<1)$$
and
\begin{align}
\label{Maz5}
&\mathbb{E}({\bf a},{\bf b},-cD_{q^{-1}})
\left\{(at,as;q)_\infty \right\} \notag \\
&\qquad = (at,as;q)_\infty\; \sum_{n=0}^\infty\; 
\frac{(a_1,a_2,\cdots,a_{\mathfrak r+1};q)_n}
{(q,b_1,b_2,\cdots,b_{\mathfrak r};q)_n}\; 
\; {}_{2}\Phi_{1}\left[
\begin{array}{rr}
q^{-n}, q/(at);\\  
\\
0;
\end{array}\,
q; as\right] \; \left(\frac{ct}{1-q}\right)^n.
\end{align} 
\end{lemma}

\begin{proof}
The first assertion (\ref{Ma2A}) of Lemma \ref{lem3.1} 
follows from (\ref{MaA}) when $s=0$. 
On the other hand, the second assertion 
(\ref{Maz5}) of Lemma \ref{lem3.1}
can be deduced from (\ref{MaA1}) by setting $\omega = 0$. 
\end{proof}

\begin{proof}[Proof of Theorems $\ref{exgend}$ and $\ref{exgend1}$] 
We suppose that the operator 
${D}_{q}$ acts upon the variable $y$. Then, in view of the formulas 
in (\ref{Propo12}), we have 
\begin{align}\label{proof9}
\sum_{n=0}^\infty \sum_{m=0}^\infty 
\phi_{n+m}^{({\bf a},{\bf b})}(x,y|q)\;
 \frac{t^n}{(q;q)_n} \; \frac{s^m}{(q;q)_m}
&=\sum_{n=0}^\infty\sum_{m=0}^\infty 
\mathbb{T}({\bf a},{\bf b},(1-q)xD_q)\left\{y^{n+m}\right\}
\; \frac{t^n}{(q;q)_n} \; \frac{s^m}{(q;q)_m}\notag \\
&= \mathbb{T}\big({\bf a},{\bf b},(1-q)xD_q\big)  
\left\{\sum_{n=0}^\infty \; 
\frac{(yt)^n}{(q;q)_n}\;\sum_{m=0}^\infty \;
\frac{(ys)^m}{(q;q)_m}\right\}\notag \\
&= \mathbb{T}\big({\bf a},{\bf b},(1-q)xD_q\big)
\left\{\frac{1}{(yt,ys;q)_\infty}\right\}.
\end{align}
The proof of the assertion (\ref{exgen}) of Theorem \ref{exgend} 
can now be completed by using the relation 
(\ref{Ma2A}) in (\ref{proof9}).
 
Similarly, we observe that
\begin{align}\label{proof10}
& \sum_{n=0}^\infty \sum_{m=0}^\infty 
(-1)^{n+m}\;q^{\binom{n}{2}+\binom{m}{2}} \;
\psi_{n+m}^{({\bf a},{\bf b})}(x,y|q)\;
\frac{t^n}{(q;q)_n}\; \frac{s^m}{(q;q)_m}\notag \\
&\qquad = 
\sum_{n=0}^\infty\sum_{m=0}^\infty 
(-1)^{n+m}\;q^{\binom{n}{2}+\binom{m}{2}} \;
\mathbb{E}\big({\bf a},{\bf b},-(1-q)xD_{q^{-1}}\big)  
\left\{y^{n+m}\right\} \; \frac{t^n}{(q;q)_n}\;
\frac{s^m}{(q;q)_m}\notag \\
&\qquad =\mathbb{E}\big({\bf a},{\bf b},-(1-q)xD_{q^{-1}}\big)  
\left\{\sum_{n=0}^\infty \;(-1)^{n}\;q^{\binom{n}{2}} 
\;\frac{(yt)^n}{(q;q)_n}\;\sum_{m=0}^\infty 
(-1)^{m}\;q^{\binom{m}{2}}\;  
\frac{(ys)^m}{(q;q)_m}\right\}\notag \\
&\qquad =\mathbb{E}\big({\bf a},{\bf b},-(1-q)xD_{q^{-1}}\big)
\left\{(yt,ys;q)_\infty \right\}, 
\end{align}
which evidently completes the proof of the assertion 
(\ref{exgen1}) of Theorem \ref{exgend1}.
\end{proof}

We next derive another Rogers type formula for   
the family of the generalized Al-Salam-Carlitz $q$-polynomials  
$\psi_{n}^{({\bf a},{\bf b})}(x,y|q)$ as follows.

\begin{theorem}{\rm $\big[$Another Rogers type formula 
for $\psi_{n}^{({\bf a}, {\bf b})}(x,y|q)\big]$}
\label{TAno} 
It is asserted that 
\begin{align} 
\label{ggen1}
&\sum_{n=0}^\infty\sum_{m=0}^\infty
(-1)^{n}\; q^{\binom{n}{2}}\;
\psi_{n+m}^{({\bf a}, 
{\bf b})}(x,y|q)\; \frac{t^n}{(q;q)_n}\; 
\frac{s^m}{(q;q)_m}\notag \\
&\qquad =\frac{(yt;q)_\infty}{(ys;q)_\infty}\;{}_{\mathfrak r+2}
\Phi_{\mathfrak r+1}\left[
\begin{array}{rr}
a_1,a_2,\cdots,a_{\mathfrak r+1},t/s;\\
\\
b_1,b_2,\cdots,b_{\mathfrak r },q/(ys);
\end{array}\,
q;  \frac{xq}{y}\right]\qquad  (|ys| <1).
\end{align}
\end{theorem}

\begin{proof} We suppose that the $q$-difference operator 
${D}_{q}$ acts upon the variable $y$. 
We then obtain 
\begin{align}\label{proof11} 
&\sum_{n=0}^\infty \sum_{m=0}^\infty (-1)^{n}\;q^{\binom{n}{2}} \;
\psi_{n+m}^{({\bf a},{\bf b})}(x,y|q) \;\frac{t^n}{(q;q)_n}\; 
\frac{s^m}{(q;q)_m}\notag \\
&\qquad = 
\sum_{n=0}^\infty\sum_{m=0}^n (-1)^{n}\; q^{\binom{n}{2}}\;
\mathbb{E}\big({\bf a},{\bf b},-(1-q)xD_{q^{-1}}\big)  
\left\{y^{n+m}\right\}\; \frac{t^n}{(q;q)_n}\;\frac{s^m}{(q;q)_m}
\notag \\
&\qquad = 
\mathbb{E}\big({\bf a},{\bf b},-(1-q)xD_{q^{-1}}\big)  
\left\{\sum_{n=0}^\infty (-1)^{n}\; q^{\binom{n}{2}}\;
\frac{(yt)^n}{(q;q)_n}\;\sum_{m=0}^\infty\;\frac{(ys)^m}{(q;q)_m} 
\right\}\notag \\
&\qquad = 
\mathbb{E}\big({\bf a},{\bf b},-(1-q)xD_{q^{-1}}\big)  
\left\{\frac{(yt;q)_\infty}{(ys;q)_\infty}\right\}.
\end{align} 
The proof of the assertion (\ref{ggen1}) of Theorem \ref{TAno} 
can now be completed by applying the formula (\ref{MaA1}) with 
$s=0$ and $\omega=s$ in (\ref{proof11}). 
\end{proof}

Another extended Rogers type formula for the family  
of the generalized Al-Salam-Carlitz $q$-polynomials 
$\Psi_{n}^{({\bf a}, {\bf b})}(x,y|q)$ is given by 
Theorem \ref{Tsxt} below.

\begin{theorem}{\rm $\big[$Another extended Rogers type 
formula for $\psi_{n}^{({\bf a},{\bf b})}(x,y|q)\big]$}
\label{Tsxt}
It is asserted that  
\begin{align}
\label{gextend}
&\sum_{n=0}^\infty \sum_{m=0}^\infty 
\sum_{k=0}^\infty 
(-1)^{n+m}\; q^{\binom{n}{2}+\binom{m}{2}}\;
\psi_{n+m+k}^{({\bf a},{\bf b})}(x,y|q)\;
\frac{t^n}{(q;q)_n}\;
\frac{s^m}{(q;q)_m}\;
\frac{\omega^k}{(q;q)_k}\notag \\
&\quad =\frac{(yt,ys;q)_\infty}{(y\omega;q)_\infty}   
\; \sum_{j=0}^\infty\; 
\frac{(a_1,a_2,\cdots,a_{\mathfrak r+1};q)_j}
{(q,b_1,b_2,\cdots,b_{\mathfrak r};q)_j}\; (xt)^j \;
{}_{3}\Phi_{2}\left[
\begin{array}{rr}
q^{-j}, q/(yt),s/\omega;\\ 
\\
q/(y\omega),0;
\end{array}\,
q; q\right] \qquad (|y\omega|<1). 
\end{align}
\end{theorem}
 
\begin{proof}  
We suppose that the operator ${D}_{q}$ acts upon 
the variable $y$. By using the formulas (\ref{Propo12}), we obtain 
\begin{align}
\label{BEL}
&\sum_{n=0}^\infty \sum_{m=0}^\infty
(-1)^{n+m}\;q^{\binom{n}{2}+\binom{m}{2}}\; 
\sum_{k=0}^\infty \psi_{n+m+k}^{({\bf a},{\bf b})}(x,y|q) 
\;\frac{t^n}{(q;q)_{n}}\;\frac{s^m}{(q;q)_{m}}\;\frac{\omega^k}{(q;q)_k}
\notag \\
&\qquad =
\sum_{n=0}^\infty\sum_{m=0}^\infty (-1)^{n+m}\;q^{\binom{n}{2}+\binom{m}{2}}\;
\sum_{k=0}^\infty \mathbb{E}\big({\bf a},{\bf b},-(1-q)xD_{q^{-1}}\big)
\left\{y^{n+m+k}\right\} 
\;\frac{t^n }{(q;q)_{n}}\;\frac{s^m}{(q;q)_{m}}\;\frac{\omega^k}{(q;q)_k}
\notag \\
& \qquad =\mathbb{E}\big({\bf a},{\bf b},-(1-q)xD_{q^{-1}}\big)
\left\{\sum_{n=0}^\infty\; (-1)^{n}\; q^{\binom{n}{2}}\; 
\frac{(yt)^n}{(q;q)_{n}}\;
\sum_{m=0}^\infty (-1)^{m}\; q^{\binom{m}{2}}\;
\frac{(ys)^m}{(q;q)_m}\;
\sum_{k=0}^\infty \frac{(y\omega)^k}{(q;q)_k}\right\}\notag \\
& \qquad =\mathbb{E}\big({\bf a},{\bf b},-(1-q)xD_{q^{-1}}\big)
\left\{\frac{(yt,ys;q)_\infty}{(y\omega;q)_\infty}\right\}.
\end{align}
Thus, in light of (\ref{MaA1}), the proof of the assertion 
(\ref{gextend}) of Theorem \ref{Tsxt} 
is completed. 
\end{proof}

\section{Srivastava-Agarwal Type Generating Functions for the 
Families of the Al-Salam-Carlitz $q$-Polynomials}

\label{qdiffff} 
In this section, we use the formulas in (\ref{Propo12}) to 
derive the Srivastava-Agarwal type generating functions involving 
the families of the Al-Salam-Carlitz $q$-polynomials 
$\phi_n^{({\bf a},{\bf b})}(x,y|q)$ and 
$\psi_n^{({\bf a},{\bf b})}(x,y|q)$. 
  
The Hahn polynomials (see \cite{Hahn049}, \cite{Hahn1949} and \cite{Hahn49})
(or, equivalently, the Al-Salam-Carlitz $q$-polynomials 
\cite{AlSalam}) are defined as follows:
\be 
\phi_n^{(a)}(x|q)=\sum_{k=0}^n \qbinomial{n}{k}{q} \;
(a;q)_k x^k \qquad \text{and} \qquad  
\psi_n^{(a)}(x|q)=\sum_{k=0}^n q^{k(k-n)}\;\qbinomial{n}{k}{q}\;
(aq^{1-k};q)_k\;x^k.
\ee 

Recently, Srivastava and Agarwal \cite{SrivastavaAgarwal} 
gave a generating function which we recall here 
as Lemma \ref{lem4.1} below.

\begin{lemma}{\rm (see \cite[Eq. (3.20)]{SrivastavaAgarwal})}
\label{lem4.1}
The following generating function holds true$:$  
\be
\label{c1sums}
\sum_{n=0}^\infty (\lambda;q)_n\; \phi_n^{(\alpha)}(x|q) 
\frac{t^n}{(q;q)_n}=\frac{(\lambda t; q)_\infty}{(t;q)_\infty}
\; {}_2\Phi_1\left[
\begin{array}{rr} 
\lambda, \alpha;\\
\\
\lambda  t; 
\end{array}\,
q; xt   
\right] \qquad (\max\{|t|, |xt|\}<1).
\ee
\end{lemma}
 
The generating function (\ref{c1sums}) is known as a 
Srivastava-Agarwal type generating function 
(see, for example, \cite{Cao-Srivastava2013}). 

In this section, we give the Srivastava-Agarwal type generating functions 
for the families of the Al-Salam-Carlitz $q$-polynomials 
$\phi_n^{({\bf a},{\bf b})}(x,y|q)$ and 
$\psi_n^{({\bf a},{\bf b})}(x,y|q)$.

\begin{theorem} 
{\rm $\big[$Srivastava-Agarwal type generating functions  
for $\phi_n^{({\bf a},{\bf b})}(x,y|q)$ and 
$\psi_n^{({\bf a},{\bf b})}(x,y|q)\big]$}
\label{TAeA21}
The following Srivastava-Agarwal type generating functions hold  
true for the families of the Al-Salam-Carlitz $q$-polynomials 
$\phi_n^{({\bf a},{\bf b})}(x,y|q)$ and 
$\psi_n^{({\bf a},{\bf b})}(x,y|q)$$:$ 

\begin{align}
\label{TAZA21}
&\sum_{n=0}^\infty (\lambda;q)_n\; \phi_n^{({\bf a},{\bf b})}(x,y|q)\;
\frac{t^n}{(q;q)_n}\notag \\
&\qquad =\frac{(\lambda yt ;q)_\infty}{(yt;q)_\infty}\;
{}_{\mathfrak r+2}\Phi_{\mathfrak r+1}\left[
\begin{array}{rr}
a_1,a_2,\cdots,a_{\mathfrak r+1},\lambda;\\
\\
b_1,b_2,\cdots,b_{\mathfrak r},\lambda yt;
\end{array}\,
q; xt\right] \qquad (|yt|<1) 
\end{align}
and
\begin{align}
\label{A21}
&\sum_{n=0}^\infty (\lambda;q)_n\;\psi_n^{({\bf a},{\bf b})}(x,y|q)
\;\frac{t^n}{(q;q)_n}\notag \\
&\qquad =\frac{(\lambda yt;q)_\infty}{(yt;q)_\infty}
\;{}_{\mathfrak r+2}\Phi_{\mathfrak r+1}\left[
\begin{array}{rr}
a_1,a_2,\cdots,a_{\mathfrak r+1},\lambda;\\
\\
b_1,b_2,\cdots,b_{\mathfrak r },q/(yt);
 \end{array}\,
q;\frac{xq}{y}\right] \qquad (|yt|<1). 
\end{align}
\end{theorem}

\begin{proof} 
We suppose that the operator  ${D}_{q}$ acts upon the variable $y$.  
According to the formulas in (\ref{Propo12}), we then obtain  
\begin{align}\label{proof15}
\sum_{n=0}^\infty (\lambda;q)_n \;
\phi_n^{({\bf a},{\bf b})}(x,y|q)\;\frac{t^n}{(q;q)_n}&=
\sum_{n=0}^\infty (\lambda;q)_n\;
\mathbb{T}\big({\bf a},{\bf b},(1-q)xD_q\big)
\left\{y^n\right\}\;\frac{t^n}{(q;q)_n}\notag \\
&=\mathbb{T}\big({\bf a}, {\bf b},(1-q)xD_q\big)
\left\{\sum_{n=0}^\infty (\lambda;q)_n \;
\frac{(yt)^n}{(q;q)_n}\right\}\notag \\
&=\mathbb{T}\big({\bf a}, {\bf b}, (1-q)xD_q\big)
\left\{\frac{(\lambda yt ;q)_\infty}{(yt;q)_\infty}\right\}.
\end{align}

Now, setting $\omega=0$ in (\ref{MaA}), we have 
\be
\mathbb{T}\big({\bf a},{\bf b},(1-q)xD_q\big)
\left\{\frac{(\lambda yt;q)_\infty}{(yt;q)_\infty}\right\}=
\frac{(\lambda yt ;q)_\infty}{(yt;q)_\infty}\;
{}_{\mathfrak r+2}\Phi_{\mathfrak r+1}\left[
\begin{array}{rr}a_1,a_2,\cdots,a_{\mathfrak r+1},\lambda;\\
\\
b_1,b_2,\cdots,b_{\mathfrak r},y\lambda t;
\end{array} \,
q; xt\right],
\ee
which, in conjunction with (\ref{proof15}), completes the proof 
of the first assertion (\ref{TAZA21}) of Theorem \ref{TAeA21}. 

The proof of the second assertion (\ref{A21}) of Theorem \ref{TAeA21} 
is much akin to that of the first assertion (\ref{TAZA21}). The details 
involved are, therefore, being omitted here. 
\end{proof}

\begin{remark}\label{rem4.1}
{\rm Upon replacing $t$ by $\lambda t,$ if we set $s=t$ in the assertion 
$(\ref{ggen1})$ of Theorem $\ref{TAno},$ we get $(\ref{A21})$.} 
\end{remark}

\begin{theorem}  {\rm $\big[$Srivastava-Agarwal type bilinear 
generating function for $\phi_n^{({\bf a},{\bf b})}(x,y|q) \big]$}
\label{TAA22}
The following Srivastava-Agarwal type bilinear generating function 
holds true for the family of the generalized 
Al-Salam-Carlitz $q$-polynomials 
$\phi_n^{({\bf a},{\bf b})}(x,y|q)$ holds true$:$
\begin{align}
\label{TAZAa1}
&\sum_{n=0}^\infty \phi_n^{({\bf a},{\bf b})}(x,y|q)
\phi_n^{(\alpha)}(\mu|q)\;\frac{t^n}{(q;q)_n}\notag \\
&\qquad =\frac{(\alpha\mu yt;q)_\infty}{(\mu yt,yt;q)_\infty}\;
\sum_{j=0}^\infty \frac{(a_1,a_2,\cdots,a_{\mathfrak r+1};q)_j}
{(q,b_1,b_2,\cdots,b_{\mathfrak r};q)_j}\;(xt)^j\;
{}_{3}\Phi_{1}\left[
\begin{array}{rr}
q^{-j}, \alpha, yt;\\
\\
\alpha\mu yt;
\end{array}\,
q;\mu q^j\right] 
\end{align}
$$(\max\{|yt|,|\mu yt|\}<1).$$
\end{theorem}

\begin{proof}
We suppose that the $q$-difference operator 
${D}_{q}$ acts upon the variable $y$. 
We then find that 
\begin{align}\label{proof17}
\sum_{n=0}^\infty \phi_n^{({\bf a},{\bf b})}(x,y|q)\; 
\phi_n^{(\alpha)}(\mu|q)\;\frac{t^n}{(q;q)_n}&=
\sum_{n=0}^\infty \mathbb{T}({\bf a},{\bf b}, 
(1-q)xD_q)\left\{y^n\right\}\phi_n^{(\alpha)}(\mu|q)\;
\frac{t^n}{(q;q)_n}\notag \\
&=\mathbb{T}\big({\bf a},{\bf b},(1-q)xD_q\big)
\left\{\sum_{n=0}^\infty \phi_n^{(\alpha)}(\mu|q)\;
\frac{(yt)^n}{(q;q)_n}\right\}\notag \\
&=\mathbb{T}\big({\bf a},{\bf b},(1-q)xD_q\big)
\left\{\frac{(\alpha\mu yt;q)_\infty}
{(yt,\mu yt;q)_\infty}\right\}.
\end{align}
The proof of the assertion (\ref{TAZAa1}) of 
Theorem \ref{TAA22} is now completed by making use  
of the relation (\ref{MaA}) in (\ref{proof17}).
\end{proof}

\begin{theorem}
{\rm $\big[$Srivastava-Agarwal type bilinear generating function 
for $\psi_n^{({\bf a},{\bf b})}(x,y|q)\big]$}
\label{Tts}
The following  Srivastava-Agarwal type bilinear generating 
function holds true for the family of the generalized  
Al-Salam-Carlitz $q$-polynomials 
$\psi_n^{({\bf a},{\bf b})}(x,y|q)$$:$
\begin{align}  
\label{ls}
&\sum_{n=0}^\infty\psi_n^{({\bf a},{\bf b})}(x,y|q)\;
\psi_n^{(\alpha)}(\mu|q)(-1)^n\;q^{\binom{n}{2}}\; 
\frac{t^n}{(q;q)_n} \notag \\
&\qquad =\frac{(\mu yt,yt;q)_\infty}{(\alpha\mu yt;q)_\infty}\;
\sum_{j=0}^\infty \frac{(a_1,a_2,\cdots,a_{\mathfrak r+1};q)_j}
{(q,b_1,b_2,\cdots, b_{\mathfrak r };q)_j}\; (xt)^j \; 
{}_{3}\Phi_{2}\left[\begin{array}{rr}
q^{-j}, q/(yt), 1/\alpha;\\ 
\\
q/(\alpha\mu yt),0;
\end{array}\,
q; q\right] \qquad (|\alpha \mu yt| <1).
\end{align} 
\end{theorem}

\begin{proof} 
We suppose that the $q$-difference operator ${D}_{q}$ 
acts upon the variable $y$.  We then obtain 
\begin{align}\label{proof18}
&\sum_{n=0}^\infty \psi_n^{({\bf a},{\bf b})}(x,y|q)\;
\psi_n^{(\alpha)}(\mu|q)\;\frac{
q^{\binom{n}{2}}\;(-t)^n}{(q;q)_n}\notag \\
&\qquad =\sum_{n=0}^\infty 
\mathbb{E}({\bf a},{\bf b},-xD_{q^{-1}})\left\{y^n\right\} \;
\psi_n^{(\alpha)}(\mu|q)\frac{q^{\binom{n}{2}}\; (-t)^n}{(q;q)_n}
\notag \\
&\qquad =\mathbb{E}\big({\bf a},{\bf b},-(1-q)xD_{q^{-1}}\big)
\left\{\sum_{n=0}^\infty \psi_n^{(\alpha)}(\mu|q)\;
\frac{(-1)^n\; q^{\frac{n}{2}}\;(yt)^n}{(q;q)_n}\right\}\notag \\
&\qquad =\mathbb{E}\big({\bf a},{\bf b},-(1-q)xD_{q^{-1}}\big)
\left\{\frac{ (yt,\mu yt;q)_\infty}{(\alpha \mu yt;q)_\infty}\right\}.
\end{align}
The proof of the assertion (\ref{ls}) of Theorem \ref{Tts} 
can now be completed by making use 
of the relation (\ref{MaA1}) in (\ref{proof18}).
\end{proof}

\section{Concluding Remarks and Observations} \label{Conclusion}

Our present investigation is motivated essentially by several recent studies
of generating functions and other results for various families of
basic (or $q$-) polynomials stemming many from the works by
Hahn (see, for example, \cite{Hahn049}, \cite{Hahn1949} and \cite{Hahn49}; 
see also Al-Salam and Carlitz \cite{AlSalam}, Srivastava and Agarwal 
\cite{SrivastavaAgarwal}, Cao and Srivastava \cite{Cao-Srivastava2013}, 
and other researchers cited herein).

In terms of the familiar $q$-difference operators $D_q$ and $D_{q^{-1}}$, 
we have first introduced two homogeneous $q$-difference   
operators $\mathbb{T}({\bf a},{\bf b},cD_q)$  and  
$\mathbb{E}({\bf a},{\bf b}, cD_{q^{-1}})$, which turn out to be 
suitable for dealing with the generalized Al-Salam-Carlitz 
$q$-polynomials $\phi_n^{({\bf a},{\bf b})}(x,y|q)$ and  
$\psi_n^{({\bf a},{\bf b})}(x,y|q)$. We have then applied each of these
two homogeneous $q$-difference operators in order to derive 
generating functions, Rogers type formulas, the 
extended Rogers type formulas and the Srivastava-Agarwal type linear 
and bilinear generating functions for each of 
these families of the generalized Al-Salam-Carlitz $q$-polynomials.

The various results, which we have presented in this paper, together
with the citations of many related earlier works are believed to 
motivate and encourage interesting further researches on the topics of
study here.

In conclusion, it should be remarked that, in a recently-published 
survey-cum-expository article, Srivastava \cite{HMS-ISTT2020} 
presented an expository overview of the classical $q$-analysis 
versus the so-called $(p,q)$-analysis
with an obviously redundant additional parameter $p$ (see, for
details, \cite[p. 340]{HMS-ISTT2020}).

\end{document}